\documentclass[12pt]{article}
\usepackage[margin=2cm]{geometry}
\usepackage{amsfonts}
\usepackage{amssymb}
\usepackage{amsthm}
\usepackage{amsmath}
\usepackage{latexsym}
\usepackage{color}
\usepackage{mathrsfs}

\begin{document}
\allowdisplaybreaks[4]
\newtheorem{lemma}{Lemma}
\newtheorem{pron}{Proposition}
\newtheorem{thm}{Theorem}
\newtheorem{Corol}{Corollary}
\newtheorem{exam}{Example}
\newtheorem{defin}{Definition}
\newtheorem{remark}{Remark}
\newtheorem{property}{Property}
\newcommand{\la}{\frac{1}{\lambda}}
\newcommand{\sectemul}{\arabic{section}}
\renewcommand{\theequation}{\sectemul.\arabic{equation}}
\renewcommand{\thepron}{\sectemul.\arabic{pron}}
\renewcommand{\thelemma}{\sectemul.\arabic{lemma}}
\renewcommand{\thethm}{\sectemul.\arabic{thm}}
\renewcommand{\theCorol}{\sectemul.\arabic{Corol}}
\renewcommand{\theexam}{\sectemul.\arabic{exam}}
\renewcommand{\thedefin}{\sectemul.\arabic{defin}}
\renewcommand{\theremark}{\sectemul.\arabic{remark}}
\renewcommand{\theproperty}{\sectemul.\arabic{property}}
\def\REF#1{\par\hangindent\parindent\indent\llap{#1\enspace}\ignorespaces}

\title{\large\bf Precise large deviations of sums of widely dependent random variables and its applications}
\author{\small Zhaolei Cui$^1$~~and~~Yuebao Wang$^2$\thanks{Research supported by National Natural Science Foundation of China
(No.s 11071182).}
\thanks{Corresponding author.
Telephone: +86 512 67422726. Fax: +86 512 65112637. E-mail:
ybwang@suda.edu.cn}
\\
{\footnotesize\it 1 School of mathematics and statistics, Changshu Institute of Technology, Suzhou 215000, China}\\
{\footnotesize\it 2 School of Mathematical Sciences, Soochow University, Suzhou 215006, China}\\
}
\date{}

\maketitle {\noindent\small {\bf Abstract }}\\

{\small  In this paper, we obtain some results on precise large deviations for non-random and random sums
of widely dependent random variables with common dominatedly varying tail distribution or consistently varying tail distribution
on $(-\infty,\infty)$.
Then we apply the results to reinsurance and insurance and give some asymptotic estimates on proportional reinsurance, random-time ruin probability
and the finite-time ruin probability.
}\\

\noindent {\small{\it Keywords:} precise large deviation; widely dependent random variables;
dominatedly varying tail distribution;  consistently varying tail distribution;
proportional reinsurance; random-time ruin probability; finite-time ruin probability}\\

\noindent {\small{2000 Mathematics Subject Classification:}
Primary 60F10; 60F05; 60G50 }\\

\section{\normalsize\bf Introduction}                 

Throughout this paper, let $X_i,\ i\geq1$ be random variables
with common distribution $F$ on $(-\infty,\infty)$, 
and $S_\tau=\sum_{k=1}^\tau X_k$ be the random sum generated by $X_i,\ i\geq1$ and a nonnegative integer valued random variable $\tau$.
In particular, if $\tau=n$, then $S_n=\sum_{k=1}^n X_k$ is called the partial sum or non-random sum for all $n\ge0$, where $S_0=0$.
Here, we call the random sum and the non-random sum together as the sum of random variables.

It is well known that, the precise large deviation of the sums of random variables is a important part of the large deviation theory.
One of the basic research objects is  the uniform asymptotic of $P(S_n>x)$ for $x\ge\gamma n$ as $n\to\infty$ where $\gamma$ is any positive constant.
Compared with the corresponding research of $\ln P(S_n>x)$, the result of $P(S_n>x)$ is more precise,
meanwhile, the latter frequently needs stronger conditions than the former for the common distribution and mutual relationship of random variables.
The purpose of this paper is to study the precise large deviation of the sums of some dependent random variables.
Thus, before giving the main results,
we first introduce the concepts of related distribution class and dependent structure of random variables in this section.

\subsection{\normalsize\bf Some distribution classes}

In this paper, all limit relations refer to $x\to\infty$ without a special statement.
Let two functions $g_1(\cdot)$ and $g_2(\cdot)$ be positive eventually.
We set $g_{1,2}=\limsup\frac{g_1(x)}{g_2(x)}$.
Then $g_1(x)=O\big(g_2(x)\big)$ means $g_{1,2}<\infty$, $g_1(x)\lesssim\ (\text{or}\ \gtrsim)\ g_2(x)$ means $g_{1,2}\ \big(\text{or}\  g_{2,1}\big)\le1$, $g_1(x)\sim g_2(x)$ means $g_{1,2}=g_{2,1}=1$, and $g_1(x)=o\big(g_2(x)\big)$ means $g_{1,2}=0$.
In addition, we denote the positive part of random variable $X$ by $X^+=X\textbf{1}_{\{X\ge0\}}$,
the negative part by $X^-=-X\textbf{1}_{\{X<0\}}$, the tail of distribution $F$ by $\overline{F}=1-F$, and the $n$-fold convolution of $F$ with itself by $F^{*n},\ n\ge1$, where $F^{*1}=F$.

We say that the distribution class
\begin{eqnarray*}
{\cal L}=\big\{F:\overline{F}(x-t)\sim\overline{F}(x)\ \text{for each}\ t\in(-\infty,\infty)\big\}
\end{eqnarray*}
is long-tailed, and the distribution class
\begin{eqnarray*}
{\cal S}=\big\{F\in\mathcal{L}:\overline{F^{*2}}(x)\sim2\overline{F}(x)\big\}.
\end{eqnarray*}
is subexponential.

In addition, we set
${\overline{F}}_{*}(y)=\liminf\frac{{\overline F}(xy)}{{\overline F}(x)},~~{\overline{F}}^{*}(y)
=\limsup\frac{{\overline F}(xy)}{{\overline F}(x)}$ for each $y>0$
and $L_F=\lim\limits_{y\downarrow1}\overline{F}_*(y)$.
Then we say that, the distribution classes
\begin{eqnarray*}
{\cal D}=\big\{F:{\overline{F}}^{*}(y)<\infty\ \text{for each}\ y>1\big\},\ \ \ {\cal C}=\{F:L_F=1\}\ \ \ \text{and}
\end{eqnarray*}
$$\mathcal{ERV}=\bigcup_{0<\alpha\le\beta<\infty}\mathcal{ERV}(\alpha,\beta)
=\bigcup_{0<\alpha\le\beta<\infty}\{F:y^{-\beta}\leq\overline{F}_*(y)\leq\overline{F}^*(y)\leq y^{-\alpha}\ \text{for each}\  y>1\}$$
are dominated varying tailed, consistently varying tailed and extended regularly varying tailed, respectively.
Particularly, if $\alpha=\beta$, then $\mathcal{ERV}(\alpha,\beta)$ reduces to the regularly varying tailed distribution class, denoted
by $\mathcal{R}_{\alpha}$.
Some properties of these distributions are introduced  by the following two propositions.

\begin{pron}\label{Proposition101}
1) The following relations  hold:\\
$$F\in\mathcal{D}\Longleftrightarrow\overline{F}_*(y)>0\ for\ each\ y>1\Longleftrightarrow\overline{F}_*(y)>0\ for\ some\ y>1
$$

2) The following inclusion relations are proper:\\
$$\mathcal {R}_{\alpha}\subset \mathcal{ERV}(\alpha,\beta)\subset\mathcal{C}\subset\mathcal{L}\cap\mathcal{D}\subset\mathcal{S}\subset\mathcal{L}.$$
\end{pron}

Further, we denote the moment index of $F$ by $I_F=\sup\{s:\int_0^\infty y^sF(dy)<\infty\}$, and
the upper Matuszewska index and lower Matuszewska index of distribution $F$ by
\begin{eqnarray*}
J_F^+=-\lim\limits_{y\to\infty}\frac{\ln \overline{F}_*(y)}{\ln y}\ \ \ \ \text{and}\ \ \ \
J_F^-=-\lim\limits_{y\to\infty}\frac{\ln \overline{F}^*(y)}{\ln y}.
\end{eqnarray*}

\begin{pron}\label{Proposition102} Let $F$ be a distribution in class $\mathcal{D}$.

1) If a constant $p<J_F^-$, then $\overline{F}(x)=o(x^{-p})$.

2) If a constant $p>J_F^+$, then $x^{-p}=o\big(\overline{F}(x)\big)$, and there exist two positive constants $C_1$ and $C_2$ such that
\begin{eqnarray*}
\frac{\overline{F}(x)}{\overline{F}(y)}\le C_1\Big(\frac{x}{y}\Big)^{p},\ \ y\ge x\ge C_2.
\end{eqnarray*}

3) $0\le J_F^-\le I_F\le J_F^+<\infty$.
\end{pron}

The above concepts and corresponding properties can be found in some references,
such as Feller \cite{F1971}, Bingham et al. \cite{BGT1987}, Cline and Samorodnitsky \cite{CS1994}, Embrechts et al. $\cite{EKM1997}$, and Tang and Tsitsiashvili $\cite{TT2003}$.

\subsection{\normalsize\bf Some dependent structures}

Based on the notion of negatively orthant dependence structure of random variables,
Wang et al. $\cite{WWG2013}$ introduced the notion of widely orthant dependence structure.

By definition, $X_i,i\ge 1$ are said to be widely upper orthant dependent (WUOD),
if for each $n\ge1$, there exists some positive number $g_U(n)$ such that,
\begin{eqnarray}\label{101}
P\Big(\bigcap^{n}_{i=1}\{X_i > x_i\}\Big)\leq g_U(n)\prod_{i=1}^nP(X_i> x_i)\ \ \ \ \ \text{for all}\ x_i\in(-\infty,\infty)\ \text{and}\ 1\le i\le n;                                      \end{eqnarray}
they are said to be widely lower orthant dependent (WLOD),
if for each $n\ge1$, there exists some  positive number $g_L(n)$ such that,
\begin{eqnarray}\label{102}
P\Big(\bigcap^{n}_{i=1}\{X_i \leq x_i\}\Big)\leq g_L(n)\prod_{i=1}^n P(X_i\leq x_i)\ \ \ \ \ \text{for all}\ x_i\in(-\infty,\infty)\ \text{and}\ 1\le i\le n;               
\end{eqnarray}
and they are said to be widely orthant dependent (WOD) if they are both WUOD and WLOD.
WUOD, WLOD and WOD structures can be called  widely dependent (WD) as a joint name.
And $g_U(n),~g_L(n)$, $n \geq 1$, are called dominating coefficients.

Clearly, $g_U(n)\ge1,~g_L(n)\ge1,n\ge2$ and $g_U(1)=g_L(1)=1$.
And some basic properties of WOD random variables are as follows, see Proposition 1.1 of Wang et al. $\cite{WWG2013}$.

\begin{pron}\label{Proposition103}
1) Let $X_i,i\geq1$ be WLOD (WUOD). If $f_i(\cdot),i\geq1$ are nondecreasing,
then $f_i(X_i),i\geq1$ are still WLOD (WUOD);
if $f_i(\cdot),i\geq1$ are nonincreasing, then $f_i(X_i),i\geq1$ are WUOD (WLOD).

2) If $X_i,i\geq1$ are nonnegative and WUOD, then for each $n\geq1$,
$$E\prod^{n}_{i=1}X_{i} \leq g_U(n)\prod^{n}_{i=1}EX_{i}.$$
In particular, for each $n\geq1$ and any $s>0$,
$$Ee^{s \sum^{n}_{i=1}X_{i}} \leq g_U(n)\prod^{n}_{i=1}Ee^{s {X_i}}.$$
\end{pron}

Further, Wang et al. $\cite{WWG2013}$ provided some examples of WD random variables,
which showed that the WD structure may contain common negatively dependent random variables,
some positively dependent random variables, and some others.
For example, when  $g_U(n)=g_L(n)=M$ for all $n\ge1$ and some positive constant $M$,
inequalities (\ref{101}) and (\ref{102}) describe extended negatively upper and lower orthant dependent (ENUOD and ENLOD) random variables, respectively.
$X_i,i\ge 1$ are said to be extended negatively orthant dependent (ENOD) if they are both ENUOD and ENLOD.
ENOD, ENUOD, ENLOD random variables are called collectively END r.v.s, see Liu $\cite{L2009}$.
More specially, if $M=1$, then we have the corresponding notions of NUOD, NLOD, NOD and ND random variables,
see, for example, Ebrahimi and Ghosh $\cite{EG1981}$ and Block et al. $\cite{BSS1982}$.

For the research on WD random variables, in addition to $\cite{WWG2013}$,
please see Wang and Cheng $\cite{WC2011}$ and Wang and Cheng $\cite{WC2020}$ for renewal theory,
Liu et al. $\cite{LGW2012}$ and Wang et al. $\cite{WCWM2012}$ for risk theory,
Wang et al. $\cite{WYL2012}$ for precise large deviation theory,
Qiu and Chen $\cite{QC2014}$, Wang et al. $\cite{WXHV2014}$, Naderi et al. \cite{NMAB2015}, Wang and Hu $\cite{WH2015}$, Chen et al. \cite{CWC2016},
Istv$\acute{a}$n et al. \cite{IPB2018} and Wu et al. \cite{WWR2019} for limit theory and statistical theory, among others.

Besides, there are also many results on ND and END random variables, which will not be described in detail here.

\subsection{\normalsize\bf Main results}

Some earlier work on the precise large deviation for non-random sums of independent and identically distributed random variables
can be found in Heyde $\cite{H1967a},\cite{H1967b},\cite{H1968}$ and Nagaev $\cite{N1969a},\cite{N1969b},\cite{N1969c}$.
With the development of the research of precise deviation, the classic case of the common distribution $F\in\mathcal{R}_{\alpha}$ on $[0,\infty)$ is attributed to Nagaev $\cite{N1973},\cite{N1979}$ ; the case of $F\in\mathcal{ERV}(\alpha,\beta)$ can be referred to Cline and Hsing $\cite{CH1991}$,
Kl\"{u}ppelberg and Mikosch $\cite{KM1997}$, Mikosch and Nagaev $\cite{MN1998}$, Tang et al. $\cite{TSJZ2001}$, among others;
and the case of $F\in\mathcal{C}$ was discussed by Ng et al. $\cite{NTYY2004}$ and Wang and Wang $\cite{WW2007}$.

For dependent random variables, there are also many corresponding results, see for example, Konstantinides and Mikosch $\cite{KM2005}$ geted some important results on the precise large deviation of sums for a stochastic recurrence equation;
later, Mikosch and Wintenberger $\cite{MW2013}$ studied the same problem for a stationary regularly varying sequence of random variables.
In this paper, we focus on another kind of dependent structure for random variables.
Theorem 3.1 of Tang $\cite{T2006}$ studies the precise large deviation of NOD random variables
with common distribution $F\in\mathcal{C}$ on $(-\infty,\infty)$.
Thereafter, Theorem 2.1 of Liu $\cite{L2009}$ replaces the condition $xF(-x)=o\big(\overline{F}(x)\big)$ in $\cite{T2006}$ with the condition
\begin{eqnarray}\label{103}
F(-x)=o\big(\overline{F}(x)\big)
\end{eqnarray}
and gets same result for ENOD random variables.
In addition, Theorem 2.1 of Wang et al. $\cite{WWC2006}$ starts to study the precise large deviation for $F\in\mathcal{D}$ on $[0,\infty)$.
Generally, Theorem 1 of Wang et al. \cite{WYL2012} discusses the case of WOD random variables with different distributions.
However, just like $\cite{T2006}$, $\cite{L2009}$ and \cite{WYL2012} still require the restriction of moment that
$E(X_1^-)^r<\infty$ for some $r>1$.
And in \cite{WYL2012}, the range of the dominating coefficients of WOD random variables is too small
to obtain the corresponding result as the following (\ref{104}).
In our view, Corollary 1 of Wang et al. \cite{WYL2012} gives such a result for identically distributed WOD random variables that,
for each integer $m\ge1$ and $v\in\big(0,m(m+1)^{-1}\big)$,
if $\max\{g_U(n),\ g_L(n)\}=O(n^{m^{-1}})$ for all $n\ge1$, then for each $\gamma>0$,
$$\overline{F}_*(v^{-1})\le\liminf\limits_{n\to\infty}\inf\limits_{x\geq\gamma n}\frac{P(S_n>x)}{n\overline{F}(x)}
\le\limsup\limits_{n\to\infty}\sup\limits_{x\geq\gamma n}\frac{P(S_n>x)}{n\overline{F}(x)}\leq\overline{F}_*^{-1}(v^{-1}).$$
Clearly, when $v\uparrow1$, $\overline{F}_*(v^{-1})\uparrow L_F$ and $m\uparrow\infty$, the latter implies that
$X_i,\ i\ge1$ be ENOD.

In the above all results, the condition that $\mu=EX_1=0$ is required.

Extending on the known results in $\cite{T2006}$, $\cite{L2009}$, $\cite{WWC2006}$ and \cite{WYL2012}, we study the precise large deviations of dependent random variables
with common distribution $F\in\mathcal{D}$ or $F\in\mathcal{C}$ on $(-\infty,\infty)$
for the two cases that $E|X_1|<\infty$ or $EX_1^-<\infty$ and $E(X_1^+)^r<\infty$ for some $r>1$.
For the former, we study it for ENOD random variables without the requirements $E(X_1^-)^r<\infty$ for some $r>1$ and $\mu=0$.
Of course, this result includes Theorem 2.1 of Liu \cite{L2009} for the case that $\mu=0$.
For the latter, we significantly expand the range of the dominating coefficients for WOD random variables with mean $\mu=0$ and give the standard and meaningful result.

\begin{thm}\label{Theorem101}
Let $X_i,i\geq1$ be ENOD random variables with common distribution $F\in\mathcal{D}$ on $(-\infty,\infty)$.
If $E|X_1|<\infty$ and the condition (\ref{103}) is satisfied, then for each $\gamma>0$,
\begin{eqnarray}\label{104}
L_F\le\liminf\limits_{n\to\infty}\inf\limits_{x\geq\gamma n}\frac{P(S_n>x)}{n\overline{F}(x)}
\le\limsup\limits_{n\to\infty}\sup\limits_{x\geq\gamma n}\frac{P(S_n>x)}{n\overline{F}(x)}\leq L_F^{-1}.\ \
\end{eqnarray}
In particular, if $F\in{\mathcal C}$, then for each $\gamma>0$,
\begin{eqnarray}\label{105}
\lim\limits_{n\to\infty}\sup\limits_{x\geq\gamma n}\Big|\frac{P(S_n>x)}{n\overline{F}(x)}-1\Big|=0.\ \
\end{eqnarray}
\end{thm}

\begin{thm}\label{Theorem102}
Let $X_i,i\geq1$ be WOD random variables with common distribution $F\in\mathcal{D}$ on $(-\infty,\infty)$.
Assume that $E(X_1^+)^r<\infty$ for some $r>1$, $\mu=0$ and (\ref{103}) is satisfied.

1) If there is some $a>0$ such that $g_U(n)=O(n^a)$, then the asymptotic upper bound of (\ref{104}) holds for each $\gamma>0$.

2) If there is some $b>0$ such that $g_L(n)=O(n^b)$ and $g_U(n)=O(n^{r-1})$, then the asymptotic lower bound of (\ref{104}) holds for each $\gamma>0$.

3) If all conditions in 1) and 2) are satisfied, then (\ref{104}) holds for each $\gamma>0$.
Further, if $F\in{\mathcal C}$, then (\ref{105}) holds for each $\gamma>0$.
\end{thm}

We prove Theorem \ref{Theorem101} and Theorem \ref{Theorem102} in Section 2.
Further, we investigate the precise large deviations for random sums of WOD random variables in Section 3.
Finally, in Section 4 we apply the above results to insurance and reinsurance
and give some asymptotic estimates on proportional reinsurance, random-time ruin probability
and the finite-time ruin probability.
In addition, we try to explain the meaning of WOD structure of random variables in a nonstandard renewal risk model.

\section{\normalsize\bf Proofs of the main results}
\setcounter{equation}{0}\setcounter{thm}{0}\setcounter{lemma}{0}

\subsection{\normalsize\bf Proof of Theorem \ref{101}}

Because $L_F=1$ for $F\in{\mathcal C}$, we just need to prove (\ref{104}) for each $\gamma>0$.
To this end, we need the following two lemmas for nonnegative random variables.

\begin{lemma}\label{Lemma201}
Let $X_i,i\geq 1$ be nonnegative WUOD random variables with common distribution $F$.
If $F\in\mathcal{D}$ and
\begin{eqnarray}\label{201}
g_U(n)n\overline{F}(n)\to0,\ \ \ \ \text{as}\ \ n\to\infty,
\end{eqnarray}
then for each $\gamma>0$,
\begin{eqnarray}\label{202}
\liminf\limits_{n\rightarrow\infty}\inf\limits_{x\geq\gamma n}\frac{P(S_n>x)}{n\overline{F}(x)}\geq1.
\end{eqnarray}
\end{lemma}

\proof According to the above conditions, it  holds uniformly for all $x\ge\gamma n$ that
\begin{eqnarray*}
&&P(S_n>x)\geq\sum_{i=1}^{n}P(X_i>x)-\sum_{1\le j<i\le n}P(X_i>x,X_j>x)\nonumber\\
&\geq&n\overline{F}(x)\big(1-g_U(n)n\overline{F}(x)\big)\nonumber\\
&\geq&n\overline{F}(x)\big(1-g_U(n)n\overline{F}(\gamma n)\big)\nonumber\\
&\sim& n\overline{F}(x),\ \ \ \ \text{as}\ \ n\to\infty.
\end{eqnarray*}
Therefore, (\ref{202}) holds.
\hfill$\square$\\

Clearly, the condition (\ref{201}) implies $x\overline{F}(x)\to0$,
which is slightly weaker than the condition that $\mu=EX_1<\infty$.

The following lemma is different from Theorem 3.1 of Ng et al. $\cite{NTYY2004}$, where the main research object is $P(S_n-n\mu>x)$.
However, the former plays an important role in the application, see Proposition \ref{Proposition401}-Proposition \ref{Proposition402} below.

\begin{lemma}\label{Lemma202}
Let $X_i,i\geq 1$ be nonnegative ENUOD random variables with common distribution $F\in{\mathcal D}$ and finite mean $\mu$.
Then for each $\gamma>0$,
\begin{eqnarray}\label{203}
1\le\liminf\limits_{n\to\infty}\inf\limits_{x\geq\gamma n}\frac{P(S_n>x)}{n\overline{F}(x)}
\le\limsup\limits_{n\to\infty}\sup\limits_{x\geq\gamma n}\frac{P(S_n>x)}{n\overline{F}(x)}\leq L_F^{-1}.\ \ \
\end{eqnarray}
In particulary, if $F\in{\mathcal C}$, then for each $\gamma>0$, (\ref{105})  holds.
\end{lemma}

\proof We also only prove (\ref{203}) for each $\gamma>0$.

First, according to Lemma \ref{Lemma201}, by $\mu<\infty$ and $g_U(n)=M$ for some $M>0$ and all $n\ge1$, we immediately give the asymptotic lower bound of (\ref{203}) for each $\gamma>0$.

Nexst, we prove the asymptotic upper bound of (\ref{203}).
There are only a few small differences between the following proof and the proof of Theorem 3.1 in Ng et al. $\cite{NTYY2004}$.
For the sake of completeness, we still give the details of the proof.
For any constant $v\in(0,1)$, we denote $\overline{X_i}=\min\{X_i,vx\},\ i\ge1$ with common distribution $F_1$,
and $\overline{S_n}=\sum\limits_{i=1}^{n}\overline{X_i},\ n\geq1$, then
\begin{eqnarray}\label{eqn7414}
P(S_n>x)\leq n\overline{F}(vx)+P(\overline{S_n}>x).
\end{eqnarray}
In order to estimate $P(\overline{S_n}>x)$ in (\ref{eqn7414}), we define an function
\begin{eqnarray*}
a(x)=\max\{-\ln\big(n\overline{F}(vx)\big),1\},\ \ \ x\in[0,\infty).
\end{eqnarray*}
By $\mu<\infty$, we know that $x\overline{F}(x)\to0$, thus
\begin{eqnarray*}
\liminf\limits_{n\to\infty}\inf\limits_{x\geq\gamma n}a(x)
\ge\lim\limits_{n\to\infty}-\ln(v\gamma)^{-1}\big(v\gamma n\overline{F}(v\gamma n)\big)=\infty.
\end{eqnarray*}
For any two constants $\tau>1$ and $\rho>\max\{J^+_F,\ \tau^{-1}\}$, we again define a function $s(\cdot)$ such that
\begin{eqnarray*}
s(x)=\frac{a(x)-\rho\tau\ln a(x)}{vx}>0,\ \ \ x\in[x_1,\infty)\ \ \text{for some}\ x_1>0.
\end{eqnarray*}
Because $\overline{X_i},i\geq1$ are also ENUOD, by Proposition \ref{Proposition102} 2) and inequality $y+1\le e^y,\ y\in[0,\infty)$,
we know that, for all integers $n$ and $x\in[x_1,\infty)$,
\begin{eqnarray}\label{eqn7415}
&&\frac{P(\overline{S_n}>x)}{n\overline{F}(vx)}\leq Me^{-s(x)x+a(x)}\big(Ee^{s(x)\overline{X_1}}\big)^n\nonumber\\
&=&Me^{-s(x)x+a(x)}\Big(\int_{0-}^{vx}(e^{s(x)y}-1)F_1(dy)+\int_{vx}^{\infty}(e^{s(x)vx}-1)F_1(dy)+1\Big)^n\nonumber\\
&\le&Me^{-s(x)x+a(x)+n\int_{0-}^{vx}(e^{s(x)y}-1)F_1(dy)}=Me^{I(n,x)}. \ \ \                    
\end{eqnarray}
Further, we estimate $n\int_{0-}^{vx}(e^{s(x)y}-1)F_1(dy)$.
According to Proposition \ref{102} 3), by inequality $e^u-1\le ue^u,\ u\in[0,\infty)$, $F\in{\cal D}$ and $\rho>J_F^+$,
we know that, there exist two positive constants $C_1$ and $C_2$ such that, for all $x\ge \max\{C_1,\ x_1\}$,
\begin{eqnarray}\label{eqn7416}
&&n\int_{0-}^{vx}(e^{s(x)y}-1)F_1(dy)=n\Big(\int_{0-}^{vxa^{-\tau}(x)}+\int_{vxa^{-\tau}(x)}^{vx}\Big)(e^{s(x)y}-1)F_1(dy)\nonumber\\
&\le&n\int_{0-}^{vxa^{-\tau}(x)}s(x)ye^{s(x)y}F(dy)+ne^{s(x)vx}\overline{F}\big(vxa^{-\tau}(x)\big)\nonumber\\
&\le&\mu ns(x) e^{s(x)vxa^{-\tau}(x)}+C_2ne^{s(x)vx}a^{\rho\tau}(x)\overline{F}(vx)\nonumber\\
&\le&\mu ns(x)e^{a^{1-\tau}(x)}+C_2
\end{eqnarray}
Then by (\ref{eqn7416}), $\tau>1$ and $v<1$, for each $\gamma>0$, we have
\begin{eqnarray}\label{eqn7417}
&&I(n,x)\le \mu ns(x)e^{a^{1-\tau}(x)}+C_2-s(x)x+a(x)\nonumber\\
&\le&\mu\gamma^{-1}v^{-1}vxs(x)e^{a^{1-\tau}(x)}+C_2-v^{-1}vxs(x)+a(x)\nonumber\\
&\le&\Big(\gamma^{-1}v^{-1}\mu e^{a^{1-\tau}(x)}-\big(v^{-1}-\frac{\rho\tau\ln a(x)}{a(x)}-1\big)\Big)a(x)+C_2\nonumber\\
&\sim&\gamma^{-1}v^{-1}\mu-(v^{-1}-1)a(x)+C_2,\ \ \ \ x\ge\gamma n\ \text{and}\ n\to\infty.
\end{eqnarray}
Thus, by (\ref{eqn7415}), (\ref{eqn7417}), $F\in{\cal D}$ and (\ref{201}), it  holds uniformly for all $x\ge\gamma n$ that
\begin{eqnarray}\label{209}
\limsup\limits_{n\to\infty}\sup\limits_{x\geq\gamma n}\frac{P(\overline{S_n}>x)}{n\overline{F}(vx)}
\le\limsup\limits_{n\to\infty}\sup\limits_{x\geq\gamma n}Me^{\gamma^{-1}v^{-1}\mu-(v^{-1}-1)a(x)+C_2}=0.
\end{eqnarray}
Therefore, by (\ref{eqn7414}), (\ref{209}) and $F\in{\cal D}$, we have
\begin{eqnarray}\label{eqn7418}
\limsup\limits_{n\to\infty}\sup\limits_{x\geq\gamma n}
\frac{P(S_n>x)}{n\overline{F}(x)}\leq\limsup\frac{\overline{F}(vx)}{\overline{F}(x)}
=\Big(\liminf\frac{\overline{F}(vxv^{-1})}{\overline{F}(vx)}\Big)^{-1}
=\overline{F}_*^{-1}(v^{-1}).
\end{eqnarray}
Then, let $v\uparrow1$ in (\ref{eqn7418}), we obtain the asymptotic upper bound in (\ref{203}) for each $\gamma>0$.
$\hspace{\fill}\Box$\\

Now we prove (\ref{104}) for each $\gamma>0$.
According to Proposition \ref{103} 1), $X_i^+,\ i\geq1$ are also ENUOD.
Further, according to Lemma \ref{Lemma202}, by $P(S_n>x)\le P(\sum_{i=1}^nX_i^+>x)$ and $EX_1^+<\infty$,
we immediately obtain the asymptotic upper bound of (\ref{104}) for each $\gamma>0$.

To prove the asymptotic lower bound of (\ref{104}) for each $\gamma>0$, we set
\begin{eqnarray}\label{210}
&A_i=\{X_i>(v+1)x,\max\limits_{1\le j\neq i\le n}X_j\le(v+1)x\},\ \ 1\le i\le n,\ n\ge1,
\end{eqnarray}
where $v$ is any constant on $(0,1)$, and deal with $P(S_n>x)$ as follows:
\begin{eqnarray}\label{eqn7420}
&&P(S_n>x)\ge P\Big(S_n>x,\bigcup_{i=1}^nA_i\Big)\nonumber\\
&=&\sum_{i=1}^nP(A_i)-\sum_{i=1}^nP(S_n\le x,A_i)=P_1(n,x)-P_2(n,x).
\end{eqnarray}

We first estimate $P_1(n,x)$. From
\begin{eqnarray*}
&&P_1(n,x)\ge n\overline{F}\big((v+1)x\big)-\sum_{i=1}^n\sum\limits_{1\le j\neq i\le n}P\big(X_l>(v+1)x,l=i,j\big)\\
&\ge&n\overline{F}\big((v+1)x\big)\big(1-Mn\overline{F}((v+1)x)\big),
\end{eqnarray*}
$x\overline{F}(x)\to0$ and $F\in\mathcal{D}$, we know that
\begin{eqnarray}\label{eqn7421}
&&\liminf\limits_{n\to\infty}\inf\limits_{x\ge\gamma n}\frac{P_1(n,x)}{n\overline{F}(x)}
=\liminf\limits_{n\to\infty}\inf\limits_{x\ge\gamma n}\frac{P_1(n,x)}{n\overline{F}\big((v+1)x\big)}\frac{\overline{F}\big((v+1)x\big)}{\overline{F}(x)}\nonumber\\
&\ge&\overline{F_*}\big((v+1)\big)\uparrow L_F,\ \ \ \ \ \text{as}\ v\downarrow0.
\end{eqnarray}

For $P_2(n,x)$, we set $\overline{X_i^-}=X_i^-\wedge vx$ with the distribution $F^-$ and $\overline{S^-_n}=\sum\limits_{i=1}^{n}\overline{X^-_i},\ \ 1\le i\le n,\ n\geq1$.
Because $\overline{X_i^-},i\ge1$ are also ENUOD
and $\overline{F^-}(x-)=o\big(\overline{F}(x)\big)$, we have
\begin{eqnarray*}
&&P_2(n,x)\le\sum_{i=1}^nP\Big(\sum_{1\le j\neq i\le n}-X_j\ge vx,A_i\Big)\nonumber\\
&\le&P\Big(\sum_{j=1}^nX_j^-\ge vx\Big)\nonumber\\
&\le&n\overline{F^-}(vx-)+P\big(\overline{S_n^-}\ge vx\big)\nonumber\\
&=&o\big(n\overline{F}(vx)\big)+Me^{-s(x)x+a(x)}E^ne^{s(x)\overline{X_1^-}}n\overline{F}(vx)\nonumber\\
&\le&o\big(n\overline{F}(vx)\big)+Me^{-s(x)x+a(x)+n\int_{0}^{vx}(e^{s(x)y}-1)F^-(dy)}n\overline{F}(vx),
\end{eqnarray*}
where the functions $a(\cdot)$ and $s(\cdot)$ are defined as Lemma \ref{Lemma202}.
Using the method similar to the proof of (\ref{eqn7416}), we have
\begin{eqnarray*}
&&n\int_{0-}^{vx}(e^{s(x)y}-1)F^-(dy)\le n\int_{0-}^{vxa^{-\tau}(x)}s(x)ye^{s(x)y}F^-(dy)\nonumber\\
&&\ \ \ \ \ \ \ \ \ \ \ \ \ \ \ \ \ \ \ +ne^{s(x)vx}\overline{F^-}\big(vxa^{-\tau}(x)\big)\nonumber\\
&\le&ns(x)e^{s(x)vxa^{-\tau}(x)}\mu_{F^-}+o\big(ne^{s(x)vx}\overline{F}\big(vxa^{-\tau}(x)\big)\nonumber\\
&\le&ns(x)e^{a^{1-\tau}(x)}\mu_{F^-}+o\big(ne^{s(x)vx}a^{\rho\tau}(x)\overline{F}(vx)\big).
\end{eqnarray*}
The rest of the proof is the same as that of Lemma \ref{Lemma202}.
Thus, for each $\gamma>0$, we have
\begin{eqnarray}\label{212}
\limsup\limits_{n\to\infty}\inf\limits_{x\ge\gamma n}\frac{P_2(n,x)}{n\overline{F}(x)}=0
\end{eqnarray}

Combining (\ref{eqn7421}) and (\ref{212}), we know that the asymptotic lower bound of (\ref{104}) holds for each $\gamma>0$.

\subsection{\normalsize\bf Proof of Theorem \ref{102}}

We use a new method to prove the theorem.
For this purpose, we first introduce some concepts and their properties.

Let $g(\cdot)$ be a positive function on $[0,\infty)$.
Further, we define some functions on $[0,\infty)$ as follows:
\begin{eqnarray}\label{eqn4110}
f_{g}(x)=\frac{1}{g(x)}\Big(\int_{0-}^x+\int_{-x-}^0\Big)|y|g(y)F(dy)=f_{g}^+(x)+f_{g}^-(x)
\end{eqnarray}
and
\begin{eqnarray}\label{eqn4111}
F_{g}(x)=\big(f_{g}^+(x)+xP(X_1>x)\big)+\big(f_{g}^-(x)+xP(X_1<-x)\big)=F_{g}^+(x)+F_{g}^-(x).
\end{eqnarray}
In addition, we say that some positive function $s(\cdot)$ almost increases to infinity, denoted by $s(x)\underline{\uparrow}\infty$,
if there are two positive constants $C$ and $x_0$ such that
$$s(x)\le C s(y)\to\infty\ \ \ \text{for all}\ x_0\le x<y.$$
Similarly, we can define concept and notation of almost decrease to zero as $s(x)\overline{\downarrow}0$.
Then, we characterize the convergence rate of the function $F_{g}^+(x)$ and the upper bound of $P(S_n>x)$.
See Proposition 1 and (2.10) in Theorem 1 of Chen et al. $\cite{CWC2016}$, respectively.

\begin{lemma}\label{Lemma203}
For $j=1,2$, let $g_j(\cdot)$ be a positive even function such that $g_j(x)\underline{\uparrow}\infty$, and let function $g(\cdot)=g_1(\cdot)g_2(\cdot)$.
If $EX_1^+g_1(X_1)<\infty$, then
\begin{eqnarray}\label{eqn4112}
&\lim g_1(x)F_{g}^+(x)=0.
\end{eqnarray}
\end{lemma}

\begin{lemma}\label{Lemma204}
Let $X_i,i\geq1$ be WUOD random variables with common distribution $F$ on $(-\infty,\infty)$ and mean $\mu=0$,
$g_1(\cdot)$ and $g_2(\cdot)$ be two even functions such that $g_j(x)\underline{\uparrow}\infty,\ i=1,2$, and $g(\cdot)=g_1(\cdot)g_2(\cdot)$.
If $EX_1^+g_1(X_1)<\infty$ and for some positive integer $m$,
\begin{eqnarray}\label{2a}
\frac{x^{m}}{g(x)}\underline{\uparrow}\infty\ \ \ \ \ \text{and}\ \ \ \ \ \ \frac{x^{m-1}}{g_1(x)}\overline{\downarrow}0,
\end{eqnarray}
then for each $\gamma>0$ and any $0<v,\theta<1$, there are two positive constants $x_0=x_0\big(F,m,v,g_1(\cdot),$ $g_2(\cdot)\big)$ and $C_0=C_0\big(F,m,v,g_1(\cdot),g_2(\cdot)\big)$ such that
\begin{eqnarray}\label{eqn4121}
P(S_n>x)\le nP(X_1>vx)+C_0g_U(n)\big(F_{g}^+(vx)\big)^{\frac{1-\theta}{v}}\ \ \ \ \ \text{for all}\ n\ge1\ \text{and}\ x\ge\gamma n\vee x_0.
\end{eqnarray}
\end{lemma}

Now, we prove Theorem \ref{102}.

1) From $E(X_1^+)^r<\infty$ for some $\gamma>1$, we can take $g_1(x)=|x|^{r-1},\ x\in(-\infty,\infty)$,
some even functions $g_2(x)\underline{\uparrow}\infty$ and positive integer $m$ satisfying (\ref{2a}).
According to Lemma \ref{Lemma203}, we have
\begin{eqnarray}\label{219}
&F_g^+(x)=o\big(g_1^{-1}(x)\big)=o\big(x^{1-r}\big).
\end{eqnarray}
And according to Lemma \ref{Lemma204}, we know that, for each $\gamma>0$ and any $0<v,\theta<1$, there are two positive constants $C_0$ and $x_0$ such that (\ref{eqn4121}) holds for all $x\ge\gamma n\vee x_0$ and all $n\ge1$.
And by Proposition \ref{Proposition102} 3), we know that $x^{-p}=o\big(\overline{F}(x)\big)$ for each $p>J_F^+$.
Further, we take $c=a\vee(r-1)$ and $v=(1-\theta)^{k+1}$ for $k=k(c,p,r,\theta)$ large enough such that
\begin{eqnarray}\label{eqn7423}
c+p-(r-1)(1-\theta)v^{-1}=c+p-(r-1)(1-\theta)^{-k}<0.
\end{eqnarray}

For all $x\ge\gamma n$, by $g_U(n)=O(n^a)$ and (\ref{eqn7423}),
it holds uniformly for all $x\ge\gamma n$ that
\begin{eqnarray}\label{eqn7424}
&&\frac{C_0g_U(n)\big(F_{g}^+(vx)\big)^{\frac{1-\theta}{v}}}{n\overline{F}(x)}=o\big(n^{a-1}x^{p-(r-1)(1-\theta)^{-k}}\big)\nonumber\\
&=&o\big(n^{-1+c+p-(r-1)(1-\theta)^{-k}}\big)\to0,\ \ \ \ \text{as}\ n\to\infty.
\end{eqnarray}
And by $F\in\mathcal{D}$ and $v=(1-\theta)^{k+1}\to1\Longleftrightarrow \theta\to0$, we know that
\begin{eqnarray}\label{eqn7425}
\limsup\limits_{n\to\infty}\sup\limits_{x\ge\gamma n}\frac{n\overline{F}(vx)}{n\overline{F}(x)}
=\Big(\liminf\frac{\overline{F}(\frac{vx}{v})}{\overline{F}(vx)}\Big)^{-1}
=F_*^{-1}(v^{-1})\to L_F^{-1},\ \ \text{as}\ v\uparrow1.
\end{eqnarray}
Combining with (\ref{eqn4121}), (\ref{eqn7424}) and (\ref{eqn7425}), we obtain the asymptotic upper bound in (\ref{104}).

2) For any $v\in(0,1)$, we set $A_i,1\le i\le n,\ n\ge1$ such as (\ref{210}).
We just deal with $P_1(n,x)$ and $P_2(n,x)$ in (\ref{eqn7420}), respectively.

For $P_1(n,x)$, by $F\in\mathcal{D},\ E(X_1^+)^r<\infty$ and $g_L(n)=O(n^{r-1})$, we have
\begin{eqnarray*}
&&P_1(n,x)\ge n\overline{F}\big((v+1)x\big)-\sum_{i=1}^n\sum\limits_{1\le j\neq i\le n}P\big(X_l>(v+1)x,l=i,j\big)\\
&\ge&n\overline{F}((v+1)x)\big(1-g_U(n)n^{-r+1}n^r\overline{F}((v+1)x)\big)\\
&\sim&n\overline{F}\big((v+1)x\big),\ \ \ \ \ \ \ x\ge\gamma n\ \text{as}\ n\to\infty.
\end{eqnarray*}
Thus, (\ref{eqn7421}) also holds for the WOD random variables.

For $P_2(n,x)$, we set $Y_i=-X_i,\ \overline{Y_i}=Y_i\wedge vx,\ i\ge1$.
By $EY_1=0$, $E\overline{Y_1}\le0$.
Thus,
\begin{eqnarray}\label{eqn7426}
&&P_2(n,x)=\sum_{i=1}^nP(S_n\le x,A_i)\le\sum_{i=1}^nP\Big(\sum_{1\le j\neq i\le n}Y_j>vx,A_i\Big)\nonumber\\
&\le&\sum_{i=1}^nP\Big(\bigcup_{1\le j\neq i\le n}\{Y_j>vx,A_i\}\Big)+\sum_{i=1}^nP\Big(\sum_{1\le j\neq i\le n}\overline{Y_j}>vx\Big)\nonumber\\
&\le&P\Big(\bigcup_{i=1}^nY_i>vx\Big)+\sum_{i=1}^nP\Big(\sum_{1\le j\neq i\le n}(\overline{Y_j}-E\overline{Y_j})>vx\Big)\nonumber\\
&=&P_{21}(n,x)+P_{22}(n,x).
\end{eqnarray}
And by $F(-x)=o\big(\overline{F}(x)\big)$ and $F\in\mathcal{D}$, it holds uniformly for all $x\ge\gamma n$ that
\begin{eqnarray}\label{eqn7427}
P_{21}(n,x)\le nP\big(Y_1>vx\big)\le nF(-vx)=o\big(n\overline{F}(x)\big),\ \ \text{as}\ n\to\infty.
\end{eqnarray}
We denote the distribution of $\overline{Y_1}-E\overline{Y_1}$ by $G$.
For each integer $n\ge1$, because $\overline{Y_j}-E\overline{Y_j},1\le j\neq i\le n$
still are WUOD random variables with the dominating coefficients $g_L(n),n\ge1$ and mean $m_G=0$.
In addition, by $E(X_1^+)^r<\infty$ and (\ref{103}), we easily know that $E\big((\overline{Y_j}-E\overline{Y_j})^+\big)^r<\infty$.
Therefore, according to Proposition \ref{Proposition103} 1), $g_L(n)=O(n^b)$, (\ref{eqn7423}) and the above proof of the asymptotic upper bound,
it holds uniformly for all $x\ge\gamma n$ that
\begin{eqnarray}\label{eqn7428}
&P_{22}(n,x)\le ng_L(n)\big(G_g^+(vx)\big)^{\frac{1-\theta}{v}}=o\big(n\overline{F}(x)\big),\ \ \text{as}\ n\to\infty.
\end{eqnarray}
Combining with (\ref{eqn7426})-(\ref{eqn7428}) and (\ref{eqn7421}), we immediately get the asymptotic lower bound in (\ref{104}).

3) The result is clearly true.

\section{\normalsize\bf Random sums}
\setcounter{equation}{0}\setcounter{thm}{0}\setcounter{lemma}{0}

For $t\ge0$, let $N(t)$ be a nonnegative integer-valued random random with finite mean $\lambda(t)=EN(t)$
and $S_{N(t)}=\sum\limits_{k=1}^{N(t)}X_k$ be the random sum generated by mutually independent $N(t)$ and $X_i,i\ge1$.
Under the following conditions:
\begin{eqnarray}\label{301}
\lambda(t)\to\infty\ \ \ \ \text{and}\ \ \ \ \ N(t)\lambda^{-1}(t)\stackrel{P}{\to}1,\ \ \ \ \text{as}\ \ t\to\infty,
\end{eqnarray}                                         
and for some $p>J_F^+$,
\begin{eqnarray}\label{302}
EN^p(t)\textbf{1}_{\{N(t)>(1+\delta)\lambda(t)\}}=o\big(\lambda(t)\big),\ \ \ \ \text{as}\ \ t\to\infty,                                \end{eqnarray}
we give the following two results, respectively.

\begin{thm}\label{Theorem301}
Let the conditions of Theorem \ref{Theorem101} be satisfied.
Further, if the conditions (\ref{301}) and (\ref{302}) are satisfied, then for each $\gamma>0$,
\begin{eqnarray}\label{303}
L_F\leq\liminf\limits_{t\to\infty}\inf\limits_{x\geq\gamma
\lambda(t)}\frac{P(S_{N(t)}>x)}{\lambda(t)\overline{F}(x)}
\leq\limsup\limits_{t\to\infty}\sup \limits_{x\geq\gamma
\lambda(t)}\frac{P(S_{N(t)}>x)}{\lambda(t)\overline{F}(x)}
\leq L_F^{-1}                                              
\end{eqnarray}
In particular, if $F\in{\mathcal C}$, then for each $\gamma>0$,
\begin{eqnarray}\label{304}
\lim\limits_{t\to\infty}\sup\limits_{x\geq\gamma \lambda(t)}\Big|\frac{P(S_{N(t)}>x)}{\lambda(t)\overline{F}(x)}-1\Big|=0.\ \
\end{eqnarray}
\end{thm}

\proof We also only need to prove (\ref{303}).
For any $0<\delta<1$, using the standard method, see, for example, Lemma 4.1-Lemma 4.3 in Kl\"{u}ppelberg and Mikosch \cite{KM1997}, we segment $P(S_{N(t)}>x)$ as follows:
\begin{eqnarray}\label{305}
&&P(S_{N(t)}>x)=\Big(\sum\limits_{0\le n<(1-\delta)\lambda(t)}+\sum\limits_{(1-\delta)\lambda(t)\leq n\leq (1+\delta)\lambda(t)}
+\sum\limits_{n>(1+\delta)\lambda(t)}\Big)P(S_n>x)P(N(t)=n)\nonumber\\
&=:&\Sigma_{11}(x,t)+\Sigma_{12}(x,t)+\Sigma_{13}(x,t).                               
\end{eqnarray}

We first estimate $\Sigma_{11}(x,t)$. By the right-hand inequality of (\ref{104}) in Theorem \ref{Theorem101} and (\ref{301}),
when $t\to\infty$, it holds uniformly for all $x\ge\gamma \lambda(t)\ge\frac{\gamma}{1-\delta}\lfloor(1-\delta)\lambda(t)\rfloor$ that
\begin{eqnarray}\label{306}
&&\Sigma_{11}(x,t)\leq P\Big(\sum_{0\le i<(1-\delta)\lambda(t)}X_i^+>x\Big)P\big(N(t)<(1-\delta)\lambda(t)\big)\nonumber\\
&\lesssim&(1-\delta)\lambda(t)\overline{F}(x)L_F^{-1}P\big(N(t)\lambda^{-1}(t)-1<-\delta)\big)\nonumber\\
&=&o\big(\lambda(t)\overline{F}(x)\big).                                 
\end{eqnarray}

Next, we deal with $\Sigma_{12}(x,t)$.
On the one hand, by the second inequality of (\ref{104}) in Theorem \ref{Theorem101},
when $t\to\infty$, it holds uniformly for all $x\ge\gamma \lambda(t)\ge\frac{\gamma}{1+\delta}n$ that
\begin{eqnarray}\label{307}
\Sigma_{12}(x,t)\lesssim\sum\limits_{(1-\delta)\lambda(t)\leq n\leq (1+\delta)\lambda(t)}n\overline{F}(x)L_F^{-1}P\big(N(t)=n\big)
\leq\lambda(t)\overline{F}(x)L_F^{-1}.                                                
\end{eqnarray}
On the other hand, by the first inequality of (\ref{104}) in Theorem \ref{Theorem101} and (\ref{301}), when $t\to\infty$,
it holds uniformly for all $x\ge\gamma \lambda(t)\ge\frac{\gamma}{1+\delta}n$ that
\begin{eqnarray}\label{308}
&&\Sigma_{12}(x,t)\gtrsim\sum\limits_{(1-\delta)\lambda(t)\leq n\leq (1+\delta)\lambda(t)}n\overline{F}(x)L_FP\big(N(t)=n\big)\nonumber\\
&\geq&(1-\delta)\lambda(t)\overline{F}(x)L_FP\big((1-\delta)\lambda(t)\leq N(t)\le(1+\delta)\lambda(t)\big)\nonumber\\
&\sim&(1-\delta)\lambda(t)\overline{F}(x)L_F\to\lambda(t)\overline{F}(x)L_F,\ \ \text{as}\ \delta\to0.                                    
\end{eqnarray}

In order to estimate $\Sigma_{13}$, we need the following lemma.

\begin{lemma}\label{Lemma301}
Let $X_i,i\geq1$ be WUOD random variables with common distribution $F$ on $(-\infty,\infty)$ such that $EX_1^+<\infty$. Then, for any $u>0$,
\begin{eqnarray}\label{309}
P(S_n>x)\leq n\overline{F}(xu^{-1})+g_U(n)(eEX_1^+ nx^{-1})^u\ \ \ \ \ \text{for all}\ x>0.            
\end{eqnarray}
\end{lemma}

\proof Let $\widetilde{X}_i=\min\{X_i^+, xu^{-1}\},\ i\geq1$ with common distribution $F_2$ on $[0,\infty)$ and $\widetilde{S}_n=\sum\limits_{i=1}^{n}\widetilde{X}_i,\ n\geq1$.
Then $\widetilde{X}_i,\ i\geq1$ are also WUOD and
\begin{eqnarray}\label{310}
P(S_n>x)\le P\Big(\sum_{i=1}^nX_i^+>x\Big)\leq n\overline{F}\big(xu^{-1}\big)+P\big(\widetilde{S}_n>x\big).         
\end{eqnarray}
Since $(e^{hy}-1)y^{-1}$ for given $h>0$ is nondecreasing and $1+y\le e^y$ for all $y>0$, we have
\begin{eqnarray}\label{311}
&&P(\widetilde{S}_n>x)\leq g_U(n)e^{-hx}\big(Ee^{h\widetilde{X}_1}\big)^n\nonumber\\
&=&g_U(n)e^{-hx}\Big(1+\int_{0-}^{xu^{-1}}(e^{hy}-1)F_2(dy)\Big)^n\nonumber\\
&\leq&g_U(n)e^{-hx}\Big(1+(e^{hxu^{-1}}-1)x^{-1}u\int_{0-}^{xu^{-1}}yF(dy)\Big)^n\nonumber\\
&\leq&g_U(n)e^{-hx}\big(1+(e^{hxu^{-1}}-1)x^{-1}uEX_1^+\big)^n\nonumber\\
&\leq&g_U(n)e^{n(e^{hxu^{-1}}-1)x^{-1}uEX_1^+-hx}.                                                 
\end{eqnarray}
Let $h=h(x,n)=ux^{-1}\ln\big(xn^{-1}(EX_1^+)^{-1}+1\big)$ in (\ref{311}), we obtain that
\begin{eqnarray}\label{312}
P(\widetilde{S}_n>x)\leq g_U(n)(eEX_1^+ nx^{-1})^u\ \ \ \ \ \ \text{for all}\ x>0.                        
\end{eqnarray}
By (\ref{310}) and (\ref{312}), we know that (\ref{309}) holds.
$\hspace{\fill}\Box$
\\

Now we continue to estimate $\Sigma_{13}$. We take $u=p>J_F^+$ in (\ref{310}), then $p>I_F\ge1$ by Proposition 1.2 3).
And according to Lemma \ref{Lemma301} and Proposition \ref{Proposition102} 2), by (\ref{302}) and $p>1$, it holds uniformly for all $x\ge\gamma \lambda(t)$ that
\begin{eqnarray}\label{313}
&&\Sigma_{13}(x,t)\leq\sum\limits_{n\geq(1+\delta)\lambda(t)}\left(n\overline{F}(xp^{-1})+M(eEX_1^+ nx^{-1})^p\right)P(N(t)=n)\nonumber\\
&=&O\big(\overline{F}(x)EN(t)\textbf{1}_{\{N(t)\geq(1+\delta)\lambda(t)\}}+x^{-p}EN^{p}(t)\textbf{1}_{\{N(t)\geq(1+\delta)\lambda(t)\}}\big)\nonumber\\
&=&o\big(\lambda(t)\overline{F}(x)\big),\ \ \ \ \ \ \text{as}\ \ \ t\to\infty.                     
\end{eqnarray}
By (\ref{305})-(\ref{308}) and (\ref{313}), we know that (\ref{303}) holds.
$\hspace{\fill}\Box$\\

Using the proof method of Theorem \ref{Theorem301}, according to Theorem \ref{Theorem102}, we can prove the following result.
Here, we omit the details of its proof.

\begin{thm}\label{Theorem302}
Let $X_i,i\geq1$ be WOD random variables with common distribution $F\in\mathcal{D}$ on $(-\infty,\infty)$.
Assume that $E(X_1^+)^r<\infty$ for some $r>1$, $\mu=0$ and (\ref{103}) is satisfied.

1) If there is some $a>0$ and $p>J_F^+$ such that $g_U(n)=O(n^a)$ and
\begin{eqnarray}\label{314}
EN^{p+a}(t)\textbf{\emph{1}}_{\{N(t)>(1+\delta)\lambda(t)\}}=o\big(\lambda(t)\big),\ \ \ \text{as}\ \ t\to\infty,
\end{eqnarray}
then the asymptotic upper bound of (\ref{303}) holds for each $\gamma>0$.

2) If there is some $b>0$ and $p>J_F^+$ such that $g_L(n)=O(n^b)$, $g_U(n)=O(n^{r-1})$ and
\begin{eqnarray}\label{315}
EN^{p+r-1}(t)\textbf{\emph{1}}_{\{N(t)>(1+\delta)\lambda(t)\}}=o\big(\lambda(t)\big),\ \ \ \text{as}\ \ t\to\infty,
\end{eqnarray}
then the asymptotic lower bound of (\ref{303}) holds for each $\gamma>0$.

3) If all conditions in 1) and 2) are satisfied, then (\ref{303}) holds for each $\gamma>0$.
Further, if $F\in{\mathcal C}$, then (\ref{304}) holds for each $\gamma>0$.
\end{thm}

At the end of this section, we point out that many counting processes satisfy the conditions (\ref{301}) and (\ref{302}),
as well as (\ref{314}) and (\ref{315}),
such as, the counting process in Lemma 2.1-Lemma 2.4 of Kl\"{u}ppelberg and Mikosch $\cite{KM1997}$,
and the quasi renewal counting process in Proposition \ref{Proposition403} and Proposition \ref{Proposition402} below.

\section{\normalsize\bf Some applications}
\setcounter{equation}{0}\setcounter{thm}{0}\setcounter{pron}{0}\setcounter{lemma}{0}\setcounter{Corol}{0}

In this section, we specially introduce a quasi renewal counting process of a nonstandard renewal risk model for some insurance companies.

In this model, let the claim sizes $Y_i, i\geq1$ be WOD random variables with the common distribution $G$ on $[0,\infty)$,
the finite mean $\mu_G$ and the dominating coefficients $g_U^G(n)$ and $g_L^G(n),\ n\ge1$;
the inter-arrival times $Z_i,i\geq1$ also be WOD random variables with common distribution $H$ on $[0,\infty)$,
the finite mean $\mu_H$ and the dominating coefficients $g_U^H(n)$ and $g_L^H(n),\ n\ge1$.
And $Y_i, i\geq1$ are independent of $Z_i,i\geq1$.
The times of successive claims, $\tau_n=\sum_{i=1}^nZ_i, n\geq1$, constitute a quasi renewal counting process
$N(t)=\sum_{n=1}^\infty \textbf{1}_{\{\tau_n\leq\  t\}}$ and the corresponding renewal function $\lambda(t)=EN(t),\ t\geq0$.
For convenience, we might as well assume that $\inf\{t:\lambda(t)>0\}=0$.
In addition, let $c>0$ be the constant interest rate and $x>0$ be the initial capital of the insurance company.
In order to ensure the normal operation of the insurance company,
a safe load condition that $\mu_G<c\mu_H$ is usually assumed.

Theorem 1.2 and Theorem 1.4 of Wang and Cheng \cite{WC2011} or Theorem 4 of Chen et al. \cite{CWC2016}
established the following results for the counting process $N(t)$:

\begin{lemma}\label{lemma4010}
In the above nonstandard renewal risk model, it is specially assumed that $Z_i,\ i\ge1$ are ENOD random variables,
then
\begin{eqnarray}\label{4010}
\lambda(t)\sim\mu_H^{-1}t\ \ \text{and}\ \ N(t)\sim\mu_H^{-1}t,\ a.s.,\ \ \text{thus}\ \ N(t)\sim\lambda(t),\ a.s.\ \ \text{as}\ \ t\to\infty.
\end{eqnarray}
\end{lemma}

In the following result, 1) is from Theorem 2.1, Theorem 2.2 and Theorem 2.4 of Wang and Cheng \cite{WC2020},
which improves the corresponding results of Theorem 1.2 and Theorem 1.4 in Wang and Cheng \cite{WC2011};
2) is from Lemma 2.2 of Wang and Cheng \cite{WC2011}.

\begin{lemma}\label{lemma401}
In the above nonstandard renewal risk model, if $EZ_1^r<\infty$ for some $r>1$,
$g_U^H(n)=o(n^a)$ for some $a>0$ and $g_L^H(n)=o(n^b)$ for some $b>0,\ n\to\infty$,
then

1) (\ref{4010})  holds;

2) (\ref{314}) and (\ref{315}) hold for any $a>0$ and $p>J_G^+$.
\end{lemma}

Here, we note that, when $N(t)=n$ for given $t>0$, because $\tau_n=\sum_{i=1}^nZ_i\le t$,
$Z_i,1\le i\le n$ are in a certain kind of negative dependence.
Therefore, the hypothesis that $Z_i,i\ge1$ are ENOD or WOD is reasonable.

In the following, we only study the case that $Y_i,i\ge1$ are ENOD random variables with common distribution $G\in\mathcal{C}$ and omit the corresponding results for WOD random variables.
It should be pointed out that the former has the non-centralization shape,
which can better describe practical problems.

\subsection{\normalsize\bf Proportional reinsurance}

In the above nonstandard renewal risk model, at time $t\ge0$, a reinsurance company covers percent $q_1\cdot 100$ of the total claim amount $S_{N(t)}^G=\sum_{i=1}^{N(t)}Y_i$ for some $q_1\in(0,1)$;
the reinsurance company benefits percent $q_2\cdot 100$ of the total income $cS_{N(t)}^H=c\sum_{i=1}^{N(t)}Z_i$ for some $q_2\in(0,1)$.
Then we respectively call the proportional reinsurance and the proportional net loss at time $t\ge0$ for the reinsurance company by
$R_{11}(t)=q_1S_{N(t)}^G$ and
$$R_{12}(t)=q_1S_{N(t)}^G-q_2cS_{N(t)}^H=q_1\sum_{i=1}^{N(t)}\big(Y_i-q_1^{-1}q_2 cZ_i\big)=:q_1\sum_{i=1}^{N(t)}X_i.$$
Let $F$ be the distribution of $X_1$ and $S_n=\sum_{i=1}^nX_i,\ n\ge1$.
Clearly, $F$ is supported on $(-\infty,\infty)$.

Now, we give some results on precise deviation of the above two objects.

\begin{pron}\label{Proposition401}
In the above nonstandard renewal risk model, let $Y_i,i\ge1$ be ENOD random variables
and $Z_i,i\ge1$ be ENOD random variables.
If $EZ_1^r<\infty$ for some $r>1$, then for each $\gamma>0$,
\begin{eqnarray}\label{401}
\lim\limits_{t\to\infty}\sup\limits_{x\geq\gamma\lambda(t)}\Big|\frac{P\big(R_{11}(t)>x\big)}{\lambda(t)\overline{G}(q_1^{-1}x)}-1\Big|=0                                            \end{eqnarray}
and
\begin{eqnarray}\label{402}
\lim\limits_{t\to\infty}\sup\limits_{x\geq\gamma\lambda(t)}\Big|\frac{P\big(R_{12}(t)>x\big)}{\lambda(t)\overline{G}(q_1^{-1}x)}-1\Big|=0.                                         \end{eqnarray}
\end{pron}

\proof We only prove (\ref{402}).
Along the same line of Theorem \ref{Theorem301}, we take any fixed constant $\delta\in(0,1)$ and split $P\big(R_{12}(t)>x\big)$ as follows:
\begin{eqnarray}\label{0403}
&&P\big(R_{12}(t)>x\big)=P\big(S_{N(t)}>q_1^{-1}x\big)\nonumber\\
&=&\Big(\sum\limits_{0\le n<(1-\delta)\lambda(t)}+\sum\limits_{(1-\delta)\lambda(t)\leq n\leq (1+\delta)\lambda(t)}
+\sum\limits_{n>(1+\delta)\lambda(t)}\Big) P\Big(\sum_{i=1}^{n}\big(Y_i-q_1^{-1}q_2 cZ_i\big)>q_1^{-1}x,N(t)=n\Big)\nonumber\\
&=:&\Sigma_{11}(x,t)+\Sigma_{12}(x,t)+\Sigma_{13}(x,t).                   
\end{eqnarray}
For $\Sigma_{11}(x,t)$, similar to the proof of (\ref{306}), according to Theorem \ref{Theorem101} and Lemma \ref{lemma4010},
it holds uniformly for all $x\ge\gamma\lambda(t)\ge\gamma q_1\lfloor(1-\delta)\lambda(t)\rfloor$ that
\begin{eqnarray}\label{0404}
\Sigma_{11}(x,t)\leq P\Big(\sum_{0\le i<(1-\delta)\lambda(t)}Y_i>q_1^{-1}x\Big)P\big(N(t)<(1-\delta)\lambda(t)\big)
=o\big(\lambda(t)\overline{G}(x)\big),\ \ \text{as}\ t\to\infty.
\end{eqnarray}
For $\Sigma_{12}(x,t)$, we denote the set
$$A_{n+1}=\Big\{\sum_{i=1}^nz_i\le t,\sum_{i=1}^{n+1}z_i>t,z_i\ge0,1\le i\le n+1\Big\},\ \ n\ge1.$$
According to (\ref{105}) of Theorem \ref{Theorem101} and dominated convergence theorem, respectively,
when $t\to\infty$, it holds uniformly for all $x\ge\gamma\lambda(t)\ge\frac{q_1\gamma}{1+\delta}n$
and $(1-\delta)\lambda(t)\leq n\leq (1+\delta)\lambda(t)$ that
\begin{eqnarray}\label{0405}
&&\Sigma_{12}(x,t)=\sum\limits_{(1-\delta)\lambda(t)\leq n\leq (1+\delta)\lambda(t)}\int_{A_{n+1}}
P\Big(\sum_{i=1}^{n}Y_i>\frac{x}{q_1}+\frac{q_2c}{q_1}\sum_{i=1}^{n}z_i\Big) P(Z_i\in dz_i,1\le i\le n+1)\nonumber\\
&\sim&\sum\limits_{(1-\delta)\lambda(t)\leq n\leq (1+\delta)\lambda(t)}\int_{A_{n+1}}
n\overline{G}\Big(\frac{x}{q_1}+\frac{q_2c}{q_1}\sum_{i=1}^{n}z_i\Big) P(Z_i\in dz_i,1\le i\le n+1)\nonumber\\
&\sim&\sum\limits_{(1-\delta)\lambda(t)\leq n\leq (1+\delta)\lambda(t)}nP\big(N(t)=n\big)\overline{G}(q_1^{-1}x) \nonumber\\
&\sim&\lambda(t)\overline{G}(q_1^{-1}x)\ \ \text{as},\ t\to\infty.                                                
\end{eqnarray}
For $\Sigma_{13}(x,t)$, by Lemma \ref{lemma401} and $EZ_1^r<\infty$ for some $r>1$,
(\ref{301}), (\ref{314}) and (\ref{315})  hold for any $a>0$ and $p>J_G^+$.
Thus, by the proof of (\ref{313}), we know that, it holds uniformly for all $x\ge\gamma \lambda(t)$ that
\begin{eqnarray}\label{0406}
\Sigma_{13}(x,t)=o\big(\lambda(t)\overline{G}(x)\big),\ \ \ \ \text{as}\ \ t\to\infty.                     
\end{eqnarray}

From (\ref{0403})-(\ref{0406}), we know that (\ref{402}) holds.
$\hspace{\fill}\Box$\\

Proposition \ref{Proposition401} generalizes and improves a corresponding result in Kl\"{u}ppelberg and Mikosch $\cite{KM1997}$.
For example, Example 5.2 of their paper gives the following conclusion:

\begin{eqnarray*}\label{404}
\lim\limits_{t\to\infty}\sup\limits_{x\geq\gamma\lambda(t)}\Big|
\frac{P\big(R_{11}(t)>x\big)}{\lambda(t)\overline{G}\big(q_1^{-1}x-\mu_G\lambda(t)\big)}-1\Big|=0
\end{eqnarray*}
for each $\gamma>q_1$ and independent identically distributed random variables $Y_i,\ i\ge1$
with common distribution $G\in\mathcal{ERV}(\alpha,\beta)$ for some pair $\beta\ge\alpha>1$.

Similarly, we can study precise deviation of some other objects, such as excess-of-loss $R_{21}(t)=\sum_{i=1}^{N(t)}(Y_i-D)^+$ and excess-of-net-loss $R_{22}(t)=\sum_{i=1}^{N(t)}\big((Y_i-D)^+-q_2 cZ_i\big)$ for some $D>0$.
Here, the research goal $R_{21}(t)$, together with $R_{11}(t)$, comes from Kl\"{u}ppelberg and Mikosch $\cite{KM1997}$.

\subsection{\normalsize\bf The random-time ruin probabilities}

In standard renewal risk model, the finite-time ruin probability at some time $t\in[0,\infty)$ is defined by
\begin{eqnarray}\label{403}
&\psi(x;t)=P\Big(\sup\limits_{0\leq n\leq N(t)}S_n>x\Big).
\end{eqnarray}
Further, let $\tau$ be a nonnegative integer valued random variable, which is independent of $Y_i,Z_i,i\ge1$.
Then the random-time ruin probability is defined by
\begin{eqnarray}\label{404}
&\psi(x;\tau)=P\Big(\sup\limits_{0\leq n\leq N(\tau)}S_n>x\Big),
\end{eqnarray}
where $q_1=q_2=1$ such that $S_n=\sum_{i=1}^n X_i=\sum_{i=1}^n (Y_i-cZ_i),\ n\ge1$.

For these concepts and related research with heavy tailed claim sizes, please refer
to Asmussen \cite{A1984}, Tang \cite{T2004}, Leipus and \v{S}iaulys \cite{LS2007}, Wang et al. \cite{WCWM2012}, and so on.

In this subsection, we use Theorem \ref{Theorem101} to give the asymptotic estimation of random-time ruin probabilities for the above nonstandard renewal risk model.

\begin{pron}\label{Proposition403}
In the above nonstandard renewal risk model, let $Y_i,i\ge1$ be ENOD random variables
and $Z_i,i\ge1$ also be ENOD random variables.
If
\begin{eqnarray}\label{408}
P\big(N(\tau)>x\big)=o\big(\overline{G}(x)\big),
\end{eqnarray}
then
\begin{eqnarray}\label{409}
P\big(\psi(x;\tau)>x\big)\sim E\lambda(\tau)\overline{G}(x).
\end{eqnarray}
\end{pron}

\proof We split $P\big(\psi(x;\tau)>x\big)$ as follows:
\begin{eqnarray}\label{410}
&&P\big(\psi(x;\tau)>x\big)=\Big(\int_{[0,\infty)\cap\{x<t\}}+\int_{[0,\infty)\cap\{x\ge t\}}\Big)P\big(\sup_{0\le n\le N(t)}S_n>x\big)
P(\tau\in dt)\nonumber\\
&=:&I_1(x)+I_2(x),\ \ \ x\ge0.
\end{eqnarray}

For $I_1(x)$, we first give the following fact.
Because $\tau$ is independent of $N(x)$ for each $x>0$, according to Lemma \ref{lemma4010}, we have
\begin{eqnarray*}
P\big(N(\tau)>2^{-1}x\big)\ge P\big(N(\tau)>2^{-1}x,\ \tau>x\big)\ge P\big(N(x)>2^{-1}x\big)P(\tau>x)\sim P(\tau>x).
\end{eqnarray*}
Then by (\ref{408}) and $G\in\mathcal{C}\subset\mathcal{D}$, we know that
\begin{eqnarray}\label{411}
I_1(x)\le P(\tau>x)\lesssim P\big(N(\tau)>2^{-1}x\big)=o\big(\overline{G}(2^{-1}x)\big)=o\big(\overline{G}(x)\big).
\end{eqnarray}

For $I_2(x)$, according to Theorem \ref{Theorem101}, by $\overline{F}(x)\sim\overline{G}(x)$, there exists an integer $n_0\ge1$ such that, when $n\ge n_0$,
\begin{eqnarray}\label{412}
\sup_{x\ge n}\max\Big\{\frac{P\big(S_n>x\big)}{n\overline{G}(x)}, \frac{P\big(S_n^G>x\big)}{n\overline{G}(x)}\Big\}\le2;
\end{eqnarray}
when $1\le n\le n_0$, by (\ref{412}),
\begin{eqnarray}\label{4120}
&&\sup_{x\ge n}\frac{P\big(S_n>x\big)}{n\overline{G}(x)}\le\sup_{x\ge n_0}\frac{P\big(S_{n}>x\big)}{n\overline{G}(x)}
+\sup_{n_0>x\ge n}\frac{P\big(S_{n}>x\big)}{n\overline{G}(x)}\nonumber\\
&\le&\sup_{x\ge n_0}n_0\frac{P\big(S_{n_0}^G>x\big)}{n_0\overline{G}(x)}+\frac{1}{\overline{G}(n_0)}\nonumber\\
&\le&2n_0\overline{G}^{-1}(n_0).
\end{eqnarray}
Combining (\ref{412}) and (\ref{4120}), we have
\begin{eqnarray}\label{41200}
\sup_{x\ge n}\max\Big\{\frac{P\big(S_n>x\big)}{n\overline{G}(x)}, \frac{P\big(S_n^G>x\big)}{n\overline{G}(x)}\Big\}
\le2n_0+\frac{1}{\overline{G}(n_0)}\ \ \ \ \text{for all}\ n\ge1.
\end{eqnarray}
Now we further split $I_2(x)$ as follows:

\begin{eqnarray}\label{413}
&&I_2(x)=\Big(\int_{\{t:\ x\ge\max\{t,n)\}}+\int_{\{t:\ 0\le t\le x<n\}}\Big)\sum_{n=1}^\infty P\big(N(t)=n\big)
P\big(S_n>x\big)P(\tau\in dt)\nonumber\\
&=:&I_{21}(x)+I_{22}(x).
\end{eqnarray}
We first deal with $I_{21}(x)$. By (\ref{408}), we have $E\lambda(\tau)<\infty$.
According to (\ref{41200}), $E\lambda(\tau)<\infty$, dominated convergence theorem and $F\in\mathcal{C}\subset\mathcal{S}$, we know that
\begin{eqnarray}\label{414}
&&\lim\frac{I_{21}(x)}{\overline{G}(x)}=\int_0^\infty\sum_{n=1}^\infty P\big(N(t)=n\big)\lim\frac{P(S_n>x)}{\overline{G}(x)}
\textbf{1}_{\{t:\ x\ge\max\{t,n)\}}(t)P(\tau\in dt)\nonumber\\
&=&\int_0^\infty\sum_{n=1}^\infty nP\big(N(t)=n\big)P(\tau\in dt)\nonumber\\
&=&E\lambda(\tau).
\end{eqnarray}
Next, we deal with $I_{22}(x)$.
By (\ref{408}), we have
\begin{eqnarray}\label{4140}
I_{22}(x)\le P\big(N(\tau)>x\big)=o\big(\overline{G}(x)\big).
\end{eqnarray}

Combining (\ref{410}), (\ref{411}), (\ref{413}), (\ref{414}) and (\ref{4140}), we know that (\ref{409})  holds.
$\hspace{\fill}\Box$\\

In particular, we take $\tau=t$ for some $t>0$ in Proposition \ref{Proposition403},
then we obtain an asymptotic formula of the finite-time ruin probability in the above nonstandard renewal risk model.

\begin{Corol}\label{Corollary401}
Under the conditions of Proposition \ref{Proposition403}, for each $t>0$,
\begin{eqnarray}\label{404}
P\big(\psi(x;t)>x\big)\sim\lambda(t)\overline{G}(x).
\end{eqnarray}
\end{Corol}

This example also shows that, $P\big(N(\tau)>x\big)=O\big(P(\tau>ax)\big)$ is not true for any $a\in(0,1)$,
and many $\tau$ satisfies the condition (\ref{408}).

\subsection{\normalsize\bf The finite-time ruin probabilities}

Because the initial capital and operating time of an insurance company are very large,
it is necessary to study the precise large deviation of the finite-time ruin probability in the above nonstandard renewal risk model.

\begin{pron}\label{Proposition402}
Under the conditions of Proposition \ref{Proposition401},
\begin{eqnarray}\label{420}
\lim\limits_{t\to\infty}\sup\limits_{x\geq\gamma\lambda(t)}\Big|\frac{\psi(x;t)}{\lambda(t)\overline{G}(x)}-1\Big|=0.                                         \end{eqnarray}
\end{pron}

\proof On the one hand, we have
\begin{eqnarray}\label{421}
&&\psi(x;t)=\sum_{n=1}^\infty P\Big(\max_{1\le m\le n}\sum_{i=1}^{m}\big(Y_i-cZ_i\big)>q_1^{-1}x,N(t)=n\Big)\nonumber\\
&\le&\sum_{n=1}^\infty P\Big(\sum_{i=1}^{n}Y_i>x\Big)P\big(N(t)=n\big)\nonumber\\
&=&P\big(S_{N(t)}^G>x\big).
\end{eqnarray}
On the other hand, we have
\begin{eqnarray}\label{422}
&&\psi(x;t)\ge\sum_{n=1}^\infty P\Big(\sum_{i=1}^{n}\big(Y_i- cZ_i\big)>x,N(t)=n\Big)=P\big(S_{N(t)}>x\big).                   
\end{eqnarray}
Then, (\ref{420}) holds after we use proof methods of Theorem \ref{301} and Proposition \ref{Proposition401} to deal with (\ref{421}) and (\ref{422}), respectively.
$\hspace{\fill}\Box$\\


\begin{thebibliography}{100}

\bibitem{A1984} Asmussen, S. (1984). Approximations for the probability of ruin within finite time. Scandinavian Actuarial Journal, 31-57.

\bibitem{BGT1987} Bingham N. H., Goldie C. M. and Teugels, J. Z. (1987). Regular Variation.
Cambridge University Press.

\bibitem{BSS1982} Block H. W., Savits T. H. and Shaked M. (1982). Some
concept of negative dependence. Ann. Probab. 10, 765-772.

\bibitem{CWC2016} Chen W., Wang Y. and Cheng D. (2016). An inequality of widely dependent random variables
and its applications. Lithuanian Mathematical Journal, 56(1), 16-31.

\bibitem{CH1991} Cline D. B. H. and Hsing T. (1991). Large deviation
probabilities for sums and maxima of random variables with heavy
or subexponential tails. Preprint, Texas A$\&$M University.

\bibitem{CS1994} Cline D. B. H. and Samorodnitsky G. (1994).
Subexponentiality of the product of independent random variables.
Stoch. Process. Appl. 49, 75-98.

\bibitem{EG1981} Ebrahimi N. and Ghosh M. (1981). Multivariate negative
dependence. Comm. Statist. Theory Meth. 10, 307-337.

\bibitem{EKM1997} Embrechts P., Kl\"{u}ppelberg C. and Mikosch T. (1997). Modelling Extremal
Events for Insurance and Finance. Springer, Berlin.

\bibitem{F1971} Feller W. An Introduction to Probability Theory and its Applications, 2nd edn. vol. 2, Wiley, New York, 1971.

\bibitem{H1967a} Heyde C. C. (1967a). A contribution to the theory of large
deviations for sums of independent random variables. Z. Wahrscheinlichkeitsth. 7, 303-308.

\bibitem{H1967b} Heyde C. C. (1967b). On large deviation problem for sums of random variables which are
not attracted to the normal law. Ann. Math. Statist. 38, 1575-1578.

\bibitem{H1968} Heyde C. C. (1968). On large deviation
probabilities in the case of attraction to a nonnormal stable law. Scankhy\={a} A 30, 253-258.

\bibitem{HCW2013} He W., Cheng D. and Wang Y. (2013). Asymptotic lower bounds
of precise large deviations with nonnegative and dependent random variables. Statist. Probab. Lett. 83, 331-338.

\bibitem{IPB2018} Istv${\rm{\acute{a}}}$n F., PecsoraBettina P. and Bettina P. (2018). General theorems on exponential and Rosenthal's inequalities
and on complete convergence. J. Math. Ineq., 12(2), 433-446.

\bibitem{JP1983} Joag-Dev K. and Proschan F. (1983). Negative association of random variables with aaapplitions. Ann. Statist. 11, 286-295.


\bibitem{KM1997} Kl\"{u}ppelberg C. and Mikosch T. (1997). Large deviations of heavy-tailed random sums with
applications in insurance and finance. J. Appl. Prob. 34, 293-308.


\bibitem{KM2005} Konstantinides D. and Mikosch T. (2005). Large deviations and ruin probabilities for solutions
to stochastic recurrence equations with heavy-tailed innovations. Ann. Probab. 33, 1992-2035.


\bibitem{L1966} Lehmann E. L. (1966). Some concepts of dependence. Ann. Math. Statist. 43, 1137-1153.

\bibitem{LS2007} Leipus, R. and \v{S}iaulys, J. (2007). Asymptotic behaviour of the finite-time ruin probability under subexponential claim sizes.
Insurance Math. Econom., 40, 498-508.

\bibitem{L2009} Liu L. (2009). Precise large deviations for dependent random variables
with heavy tails. Statist. Probab. Lett. 79, 1290-1298.

\bibitem{LGW2012} Liu X., Gao Q. and Wang Y. (2012). A note on a dependent
risk model with constant interest rate. Statist. Probab. Lett. 82, 707-712.

\bibitem{M1992} Matula P. (1992). A note on the almost sure convergence of sums of negatively
dependent random variables. Statist. Prob. Lett. 15, 209-213.

\bibitem{MW2013} Mikosch T. and Wintenberger O. (2013). Precise large deviations for dependent regularly varying sequences.
Probability Theory and Related Fields volume 156, 851-887.

\bibitem{MN1998} Mikosch T. and Nagaev A.V. (1998). Large deviations of heavy-tailed sums with applications in insurance. Extremes 1, 81-110.

\bibitem{NMAB2015} Naderi H., Matula P., Amini M. and Bozorgnia A. (2015). On stochastic dominance and the strong law of large
numbers for dependent random variables. RACSAM, 110(2), 771-782.

\bibitem{N1969a} Nagaev A. V. (1969a). Integral limit theorems for large deviations when Cram\'{e}r's condition is not
fulfilled.I. Theory Prob. Appl. 14, 51-64.

\bibitem{N1969b} Nagaev A. V.(1969b). Integral limit theorems for large deviations when
Cram\'{e}r's condition is not fulfilled.II. Theory Prob. Appl. 14, 193-208.

\bibitem{N1969c} Nagaev A. V. (1969c). Limit theorems for
large deviations when Cram\'{e}r's conditions are violated.
Fiz-Mat. Nauk. 7, 17-22 (in Russian).

\bibitem{N1973} Nagaev A. V. (1973). Large deviations for sums of independent random variables.
In Trans. Sixth Prague Conf. Inf. Theory Statist. Decision Functions Random Process, Academia, 657-674.

\bibitem{N1979} Nagaev A. V. (1979). Large deviations of sums of independent random variables. Ann. Prob. 7, 745-789.

\bibitem{NTYY2004} Ng K. W., Tang Q., Yan J. and Yang H. (2004). Precise large
deviations for sums of random variables with consistently varying tails. J. Appl. Prob. 41, 93-107.

\bibitem{QC2014} Qiu D. and Chen P. (2014). Complete moment convergence for weighted sums of widely orthant dependent random variables.
Acta Math. Sin. (Engl. Ser.), 30, 1539-1548.



\bibitem{T2004} Tang Q. (2004). Asymptotics for the finite time ruin probability in the renewal model with consistent variation.
Stochastic Models, 20, 281-297.


\bibitem{T2006} Tang Q. (2006). Insensitivity to negative dependence of the asymptotic behavior of precise deviations.
Electron Journal of Probabilityb, 11, 107-120.

\bibitem{TSJZ2001} Tang Q. Su C. Jiang T. and Zhang J.S. (2001). Large deviations for
heavy-tailed random sums in compound renewal model. Statist. Prob.
Lett. 52, 91-100.

\bibitem{TT2003} Tang Q. and Tsitsiashvili G. (2003). Precise estimates for the ruin probability in finite horizon in a discrete-time model with
heavy-tailed insurance and financial risks. Stochastic Process and its Application, 108, 299-325.

\bibitem{WC2011} Wang Y. and Cheng D. (2011). Basic renewal theorems for a random walk with
widely dependent increments and their applications. J. Math. Anal. Appl. 384, 597-606.

\bibitem{WC2020} Wang Y. and Cheng D. (2020). Elementary renewal theorems for widely dependent random variables
with applications to precise large deviations. Communications in Statistics-Theory and Methods, 49 (14), 3352-3374.

\bibitem{WCWM2012} Wang Y., Cui Z., Wang K. and Ma X. (2012). Uniform asymptotics
of the finite-time ruin probability for all times. J. Math. Anal. Appl. 390, 208-223.

\bibitem{WH2015} Wang X. and Hu S. (2015). The consistency of the nearest neighbor estimator
of the density function based on WOD samples. J. Math. Anal. Appl. 429, 497-512.

\bibitem{WW2007} Wang S. and Wang W. (2007). Precise large deviations for sums of
random variables with consistently varying tails in multi-risk models.
J. Appl. Probab., 44, 4, 889-900.

\bibitem{WWC2006} Wang Y., Wang K. and Cheng D. (2006). Precise large deviations for sums of negatively associated random variables
with common dominatedly varying tails. Acta Mathematica Sinica, English Series, 22(6), 1725-1734.

\bibitem{WWG2013} Wang K., Wang Y. and Gao Q. (2013). Uniform asymptotics
for the finite-time ruin probability of a dependent risk model
with a constant interest rate. Methodol. Comput. Appl. Probab. 15, 109-124.

\bibitem{WXHV2014} Wang X., Xu C., Hu T., Volodin A. and Hu S. (2014). On complete convergence
for widely orthant-dependent random variables and its applications in nonparametric
regression models. Test, 23, 607-629.

\bibitem{WYL2012} Wang K., Yang Y. and Lin J. (2012). Precise large deviations for widely orthant dependent random variables
with dominatedly varying tails. Frontiers of Mathematics in China, 7, 919-932.


\bibitem{WWR2019} Wu Y., Wang X. and Rosalsky A. (2019). Complete Moment Convergence for the Dependent Linear Processes with Random Coefficients
Acta Mathematica Sinica, English Series, 35(8), 1321-1333.


\end{thebibliography}
\end{document}